\documentclass[a4paper,draft,reqno,12pt]{amsart}
\usepackage[english]{babel}
\usepackage{amsmath}
\usepackage{amssymb}
\usepackage{amscd}
\usepackage{amsthm}
\usepackage{euscript}
\usepackage{tikz}
\usepackage{comment} 
\newtheorem{proposition}{Proposition}

\newtheorem{theorem}{Theorem}
\newtheorem*{theorem*}{Theorem}

\theoremstyle{definition}
\newtheorem{definition}{Definition}

\theoremstyle{remark}

\DeclareMathOperator{\Hom}{Hom}
\DeclareMathOperator{\Aut}{Aut}

\DeclareMathOperator{\GL}{GL}

\def\Ker{{\rm Ker}\,}

\def\BG{{\mathbb G}}
\def\BK{{\mathbb K}}

\def\BG{{\mathbb G}}

\def\BK{{\mathbb K}}

\def\BT{{\mathbb T}}
\def\BZ{{\mathbb Z}}

\def\BQ{{\mathbb Q}}

\def\CC{\mathcal{C}}
\def\CE{\mathcal{E}}

\def\CV{\mathcal{V}}

\def\Inn{\mathrm{Inn}}

\def\MFR{\mathfrak{R}}

\title[Automorphism Groups of Reductive and Root Monoids]{Automorphism Groups of Reductive and Root Monoids}

\thanks{The paper was supported by the grant RSF 25-71-00070.}

\author{Anton Shafarevich}
\email{shafarevich.a@gmail.com}
\address{
Lomonosov Moscow State University, Faculty of Mechanics and Mathematics, Department of Higher Algebra, Leninskie Gory 1, Moscow, 119991 Russia;
\linebreak
and
\linebreak
HSE University, Faculty of Computer Science, Pokrovsky Boulevard 11, Moscow, 109028, Russia}

\begin{document}
\maketitle

\begin{abstract}

We describe the automorphism groups of reductive monoids and of root monoids with active groups of invertible elements introduced in \cite{NZ}.

\end{abstract}

\section{Introduction}
 An algebraic monoid is an algebraic variety $X$ equipped with an associative multiplication
$$
X\times X\to X
$$
given by a morphism of varieties and possessing an identity element. Although algebraic monoids have been studied extensively (see, for example, \cite{Br1,Br2,Pu,Re,Ri1}), comparatively little attention has been paid to their automorphism groups. One of results can be found in \cite{Pu2}.

If $X$ is an alebraic monoid then the group of invertible elements $G(X)$ is an open subset in $X$. If $G(X)$ is a linear group then $X$ is an affine variety, see \cite{Ri2}. We call such monoids linear. When $G(X)$ is a reductive group the monoid $X$ is called reductive. Reductive monoids were classified in~\cite{Ri1,Vin}. 

In recent years, substantial progress has been made in the classification of non-reductive linear algebraic monoids. In \cite{Bi,DZ}, all two-dimensional linear algebraic monoids were classified. In \cite{AAZ} and \cite{ABZ}, all monoid structures on the affine three-dimensional space were described. Finally, \cite{Za} contains a classification of algebraic monoids whose unit group is a solvable algebraic group of corank~$1$.

Any connected solvable algebraic group $G$ is isomorphic to a semidirect product $
G=T\ltimes_\varphi U,$
where $T$ is an algebraic torus, $U$ is a unipotent group, and
$
\varphi:T\to \Aut(U)
$
is a homomorphism. In \cite{Za}, the notion of an \emph{active} solvable algebraic group was introduced. Namely, the group $G$ is called active if
$
\dim \operatorname{Im}\varphi=\dim U.
$
If $G$ is active, then every algebraic monoid $X$ with unit group $G(X)=G$ is a toric variety.

In \cite{NZ}, the authors constructed a procedure that endows  affine toric varieties with the structure of an algebraic monoid. The monoids obtained in this way are called \emph{root monoids}. The same paper proves that every algebraic monoid with an active group of invertible elements is a root monoid.

This paper is devoted to the study of automorphism groups of linear algebraic monoids. In Section~\ref{Sec2}, we describe the automorphism groups of reductive monoids. Our description is obtained naturally from the known structure of automorphism groups of connected reductive algebraic groups.

In Section~\ref{Sec3}, we recall the construction of root monoids and describe the automorphism groups of root monoids with active groups of invertible elements. 

Throughout the paper, we assume that the ground field $\BK$ is algebraically closed and of characteristic zero. We denote by $\BG_a$ the additive group of the field $\BK$, and by $\BG_m$ its multiplicative group.

\section{Reductive monoids}\label{Sec2}

In this section we describe the automorphism groups of reductive algebraic monoids. Our description is based on the structure of automorphism groups of connected reductive algebraic groups, which we briefly recall here. Details can be found in \cite[Section~23]{Mil}.

Let $G$ be a connected reductive group, $T$ its radical, and $G_s=[G,G]$ its derived subgroup, which is semisimple. Then $G$ is an almost direct product of $T$ and $G_s$. Any automorphism of $G$ preserves both $T$ and $G_s$. The intersection $Z=T\cap G_s$ is a finite central subgroup of~$G$.

Let $k=\dim T$. Since $T\simeq (\BG_m)^k$, its automorphism group is naturally isomorphic to $\GL_k(\BZ)$. Choosing coordinates $t_1,\dots,t_k$ on $T$, a matrix
$$
A=(a_{ij})\in \GL_k(\BZ)
$$
acts on $T$ by
$$
A(t_1,\dots,t_k)
=
\left(
t_1^{a_{11}}\cdots t_k^{a_{1k}},
\dots,
t_1^{a_{k1}}\cdots t_k^{a_{kk}}
\right).
$$

We now describe the automorphism group of the semisimple factor $G_s$. Choose a Borel subgroup $B_s\subseteq G_s$ and a maximal torus $T_s\subseteq B_s$. Denote by
$$
X^*(T_s)=\Hom(T_s,\BG_m)
$$
the character lattice of $T_s$. The choice of $T_s$ determines the root system
$$
\Delta\subseteq X^*(T_s),
$$
while the choice of $B_s$ determines a system of positive roots, and hence a system of simple roots
$$
\Pi=\{\alpha_1,\ldots,\alpha_r\}\subseteq \Delta.
$$
We denote by $\mathcal D$ the corresponding Dynkin diagram. Note that the root system $\Delta$ may be reducible, and therefore the Dynkin diagram $\mathcal D$ may be disconnected.

Let $\mathfrak g_s$ and $\mathfrak t_s$ be the Lie algebras of $G_s$ and $T_s$, respectively. The lattice $X^*(T_s)$ embeds naturally into $\mathfrak t_s^*$. Every automorphism of the Dynkin diagram $\mathcal D$ induces an automorphism of the root system $\Delta$, hence an automorphism of the Lie algebra $\mathfrak g_s$ preserving $\mathfrak t_s$. Such an automorphism integrates to an automorphism of the group $G_s$ if and only if it preserves the lattice $X^*(T_s)$. We denote by
$$
\Aut(G_s,\mathcal D)
$$
the subgroup of diagram automorphisms satisfying this condition. Thus $\Aut(G_s,\mathcal D)$ may be regarded as a subgroup of $\Aut(G_s)$.

The group $\Aut(G_s)$ contains the normal subgroup $\Inn(G_s)$ of inner automorphisms, and
$$
\Inn(G_s)\simeq G_s/Z(G_s),
$$
where $Z(G_s)$ denotes the center of $G_s$. Moreover,
$$
\Aut(G_s)
\simeq
\Inn(G_s)\rtimes \Aut(G_s,\mathcal D).
$$

Let $\BT=T\cdot T_s$ be a maximal torus of $G$. Its character lattice $X^*(\BT)$ can be identified with the following sublattice of $X^*(T)\oplus X^*(T_s)$:
$$
X^*(\BT)
=
\{
(\chi_1,\chi_2)\in X^*(T)\oplus X^*(T_s)
\mid
\chi_1(z)=\chi_2(z),\ \forall z\in Z
\}.
$$
In particular, $X^*(\BT)$ is a sublattice of finite index. The set of roots $\Delta$ naturally embeds into $X^*(\BT)$.

The group $\GL_k(\BZ)\times \Aut(G_s,\mathcal D)$ acts naturally on $X^*(T)\oplus X^*(T_s)$. This defines a subgroup
$$
\Aut(G,\mathcal D)
=
\{
(\varphi,\psi)\in \GL_k(\BZ)\times \Aut(G_s,\mathcal D)
\mid
(\varphi,\psi)(X^*(\BT))=X^*(\BT)
\}.
$$
We therefore obtain
$$
\Aut(G)
\simeq
\Inn(G)\rtimes \Aut(G,\mathcal D)
\simeq
G_s/Z(G_s)\rtimes \Aut(G,\mathcal D).
$$

Let
$$
X_*(\BT)=\Hom(\BG_m,\BT)
$$
be the lattice of cocharacters of $\BT$, and set
$$
\CE=X_*(\BT)\otimes_{\BZ}\BQ,
\qquad
\CE^*=X^*(\BT)\otimes_{\BZ}\BQ.
$$
There is a natural pairing
$$
\langle\ ,\ \rangle:\CE^*\times \CE\to \BQ.
$$

Let $\alpha_1^\vee,\ldots,\alpha_r^\vee\in \CE$ be the simple coroots. Define
$$
\Lambda_+
=
\{
\lambda\in X^*(\BT)
\mid
\langle\lambda,\alpha_i^\vee\rangle\geq 0,\ i=1,\ldots,r
\}
\subseteq \CE^*
$$
to be the set of dominant weights. The closure of the negative Weyl chamber is
$$
\CV
=
\{
v\in \CE
\mid
\langle\alpha_i,v\rangle\leq 0,\ i=1,\ldots,r
\}
\subseteq \CE.
$$

Let $X$ be an algebraic monoid with unit group $G$. Then $G$ is embedded into $X$ as an open $G\times G$-stable subset. Consequently, $\BK[X]$ is a $G\times G$-invariant subalgebra of $\BK[G]$.

As a $G\times G$-module, the algebra $\BK[G]$ decomposes as
$$
\BK[G]
=
\bigoplus_{\lambda\in \Lambda_+}
V(\lambda)\otimes V(\lambda)^*,
$$
where $V(\lambda)$ denotes the irreducible $G$-module of highest weight $\lambda$.

\begin{theorem}[{\cite[Theorem~4]{Ri1}, \cite[Theorem~2]{Vin}}]
There is a one-to-one correspondence between normal affine algebraic monoids with unit group $G$ and strictly convex cones $\CC$ in $\CE$ generated by $\alpha_1^\vee,\ldots,\alpha_r^\vee$ together with a finite subset of $\CV$.

Under this correspondence, a cone $\CC$ determines a monoid $X$ whose algebra of regular functions is
$$
\BK[X]
=
\bigoplus_{\lambda\in \CC^\vee\cap X^*(\BT)}
V(\lambda)\otimes V(\lambda)^*
\subseteq
\BK[G],
$$
where $\CC^\vee\subseteq \CE^*$ is the dual cone to $\CC$.
\end{theorem}

Any automorphism of the monoid $X$ preserves the group of invertible elements $G$. Thus $\Aut(X)$ may be regarded as a subgroup of $\Aut(G)$.

For every element $g\in G$, the inner automorphism
$$
h\mapsto ghg^{-1},
\qquad
h\in G,
$$
extends naturally to an automorphism
$$
x\mapsto gxg^{-1},
\qquad
x\in X,
$$
of the monoid $X$. Consequently, the group $\Inn(G)$ embeds naturally into $\Aut(X)$.

Let $\rho\in \Aut(G)$ and let
$$
\rho^*:\BK[G]\to \BK[G]
$$
be the induced automorphism of the algebra of regular functions. Then $\rho$ extends to an automorphism of $X$ if and only if $\rho^*$ preserves the subalgebra
$$
\BK[X]\subseteq \BK[G].
$$

Suppose now that $\rho\in \Aut(G,\mathcal D)$. Then $\rho$ induces an automorphism of the lattice $X^*(\BT)$, which we denote by the same letter. Since
$$
\BK[X]
=
\bigoplus_{\lambda\in \CC^\vee\cap X^*(\BT)}
V(\lambda)\otimes V(\lambda)^*,
$$
the automorphism $\rho^*$ preserves $\BK[X]$ if and only if $\rho$ preserves $\CC^\vee$.

Define
$$
\Aut(G,\mathcal D,\CC)
=
\{
\rho\in \Aut(G,\mathcal D)
\mid
\rho(\CC^\vee)=\CC^\vee
\}.
$$

We obtain the following theorem.

\begin{theorem}
Let $\CC$ be a cone in $\CE$ generated by $\alpha_1^\vee,\ldots,\alpha_r^\vee$ together with a finite subset of $\CV$, and let $X$ be the corresponding algebraic monoid. Then
$$
\Aut(X)
\simeq
\Inn(G)\rtimes \Aut(G,\mathcal D,\CC)
\simeq
G_s/Z(G_s)\rtimes \Aut(G,\mathcal D,\CC).
$$
\end{theorem}

\section{Root monoids}\label{Sec3}
In this section we recall the definition of a root monoid introduced in \cite{NZ} and describe the automorphism groups of active root monoids.

Throughout this section, we fix an algebraic torus $\BT$ of dimension $n$. Following the notation of \cite{CLS}, we denote the character and cocharacter lattices of $\BT$ by $M$ and $N$, respectively. We also write
$$
M_\BQ=M\otimes_\BZ \BQ,
\qquad
N_\BQ=N\otimes_\BZ \BQ
$$
for the corresponding rational vector spaces. For each lattice point $m\in M$, we denote by $\chi^m$ the corresponding character of~$\BT$.

Let $\sigma \subseteq N_\BQ$ be a strongly convex rational polyhedral cone, and let $\sigma^{\vee} \subseteq M_\BQ$ be its dual cone. We denote by $X_\sigma$ the corresponding affine toric variety. If
$$
S_\sigma=M\cap \sigma^\vee,
$$
then
$$
\BK[X_\sigma]
=
\bigoplus_{m\in S_\sigma}\BK\chi^m
\subseteq
\BK[\BT].
$$

For a face $\tau$ of $\sigma$, we denote by $O_\tau$ the corresponding $\BT$-orbit in $X_\sigma$. The affine toric variety $X_\tau$ embeds naturally into $X_\sigma$ as the minimal open $\BT$-invariant affine chart containing~$O_\tau$.

A face $\tau$ is called regular if the primitive vectors along the rays of $\tau$ form part of a basis of the lattice $N$.

\begin{definition}
Let $p_1,\ldots,p_l\in N$ be the primitive vectors generating the rays of $\sigma$. For each $i$, define
$$
\MFR_i
=
\{
e\in M
\mid
\langle e,p_i\rangle=-1,
\ \langle e,p_j\rangle\geq 0
\ \text{for all}\ j\neq i
\}.
$$
The elements of the set
$$
\MFR=\bigsqcup_{i=1}^l \MFR_i
$$
are called \emph{Demazure roots}.
\end{definition}

Each Demazure root $e\in \MFR_i$ defines a derivation
$$
\delta_e:\BK[X_\sigma]\to \BK[X_\sigma]
$$
given on characters by
$$
\delta_e(\chi^m)
=
\langle m,p_i\rangle \chi^{m+e},
\qquad
m\in S_\sigma.
$$
This derivation is locally nilpotent, that is, for every $f\in \BK[X_\sigma]$ there exists an integer $q\geq 0$ such that $
\delta_e^q(f)=0.
$ The derivation $\delta_e$ defines an action of the additive group $\BG_a$ on~$X_\sigma$. The induced action on $\BK[X_\sigma]$ is given by
$$
s\cdot f
=
\exp(s\delta_e)(f)
=
\sum_{i\geq 0}
\frac{s^i\delta_e^i(f)}{i!},
\qquad
s\in \BG_a,\ f\in \BK[X_\sigma].
$$

\begin{definition}
Let $\tau$ be a $k$-dimensional regular face of $\sigma$, and let $p_1,\ldots,p_k$ be the primitive vectors generating the rays of $\tau$. A collection of Demazure roots
$$
\{e_1^{(r)},e_2^{(r)} \mid r=1,\ldots,k\}
$$
is called \emph{compatible} with $\tau$ if
$$
\langle p_i,e_1^{(r)}\rangle
=
\langle p_i,e_2^{(r)}\rangle
=
-\delta_{ri},
$$
for all $r,i=1,\ldots,k$, where $\delta_{ri}$ denotes the Kronecker delta.
\end{definition}

For a regular face $\tau$, define
$$
\tau^\perp
=
\{
m\in M_\BQ
\mid
\langle m,v\rangle =0,\ \forall v\in \tau
\}.
$$
Set $
M(\tau)=M\cap \tau^\perp.$ We define an algebraic torus
$$
T=\Hom(M(\tau),\BG_m).
$$
This torus can be naturally identified with the orbit $O_\tau$. Let $
\{e_1^{(r)},e_2^{(r)}\mid r=1,\ldots,k\}
$
be a collection of Demazure roots compatible with $\tau$. Then
$e_1^{(r)}-e_2^{(r)}\in \tau^\perp,\ 
r=1,\ldots,k.
$ Hence we obtain characters
$
\chi_r=e_2^{(r)}-e_1^{(r)}
$
of the torus $T$.

The collection
$$
\overline{\chi}=\{\chi_1,\ldots,\chi_k\}
$$
defines a homomorphism
$$
\overline{\chi}:T\to \Aut(\BG_a^k)\simeq \GL_k(\BK),
$$
given by
$$
t\mapsto
\begin{pmatrix}
\chi_1(t) & 0 & \cdots & 0 \\
0 & \chi_2(t) & \cdots & 0 \\
\vdots & \vdots & \ddots & \vdots \\
0 & 0 & \cdots & \chi_k(t)
\end{pmatrix}.
$$

We thus obtain a semidirect product
$$
G_{\overline{\chi}}
=
\BG_a^k\rtimes_{\overline{\chi}} T.
$$
Its multiplication is given by
$$
(\alpha_1,\ldots,\alpha_k,t)\cdot
(\alpha_1',\ldots,\alpha_k',t')
=
(
\alpha_1+\chi_1(t)\alpha_1',
\ldots,
\alpha_k+\chi_k(t)\alpha_k',
tt'
),
$$
where
$$
(\alpha_1,\ldots,\alpha_k),
(\alpha_1',\ldots,\alpha_k')
\in \BG_a^k,
\qquad
t,t'\in T.
$$

\begin{theorem}[{\cite[Theorem~4.2]{NZ}}]\label{ThNZ}
Let $\tau$ be a regular face of $\sigma$, and let
$$
E = \{e_1^{(r)},e_2^{(r)} \mid r=1,\ldots,k\}
$$
be a collection of Demazure roots compatible with $\tau$. Then the map
$$
\BK[X_\sigma]
\to
\BK[X_\sigma]\otimes \BK[X_\sigma]
$$
given by
$$
\chi^u
\mapsto
\chi^u\otimes \chi^u
\prod_{r=1}^k
\left(
1\otimes \chi^{e_1^{(r)}}
+
\chi^{e_2^{(r)}}\otimes 1
\right)^{\langle p_r,u\rangle}
$$
defines a comultiplication on $X_\sigma$ whose group of invertible elements is $
X_\tau\simeq G_{\overline{\chi}}.$
\end{theorem}

\begin{definition}[{\cite[Definition~4.3]{NZ}}]
The monoids arising in Theorem~\ref{ThNZ} are called \emph{root monoids}. We denote such a monoid by $X_{\sigma,E}$.
\end{definition}

Let us describe the isomorphism $
G_{\overline{\chi}}\simeq X_\tau$
from Theorem~\ref{ThNZ}. The vectors $p_1,\ldots,p_k$ can be extended to a basis
$$
p_1,\ldots,p_k,p'_{k+1},\ldots,p'_n
$$
of the lattice $N$. Let $
q_1,\ldots,q_n
$
be the dual basis in $M$. Then
$$
\BK[X_\tau]
=
\BK[
\chi^{-e_1^{(1)}},
\ldots,
\chi^{-e_1^{(k)}},
\chi^{\pm q_{k+1}},
\ldots,
\chi^{\pm q_n}
]
$$
and similarly
$$
\BK[X_\tau]
=
\BK[
\chi^{-e_2^{(1)}},
\ldots,
\chi^{-e_2^{(k)}},
\chi^{\pm q_{k+1}},
\ldots,
\chi^{\pm q_n}
].
$$

The elements $q_{k+1},\ldots,q_n$ form a basis of the lattice $M(\tau)$. Let
$$
t_1,\ldots,t_{n-k}
$$
be the corresponding coordinates on the torus $T$. Then there is an isomorphism
$$
\BK[G_{\overline{\chi}}]
\simeq
\BK[X_\tau]
$$
sending $\alpha_i$ to $\chi^{-e_1^{(i)}}$ and $t_i$ to $\chi^{q_{k+i}}$. The corresponding dual morphism defines the isomorphism $
G_{\overline{\chi}}\simeq X_\tau.
$
\begin{definition}[{\cite[Definition~5.1]{NZ}}]
Let $
H=S\ltimes_\theta U
$
be a solvable algebraic group, where $S$ is a torus, $U$ is a unipotent group, and
$$
\theta:S\to \Aut(U)
$$
is a homomorphism. The group $H$ is called \emph{active} if
$
\dim \operatorname{Im}\theta=\dim U.
$
\end{definition}

\begin{theorem}[{\cite[Theorem~5.5]{NZ}}]
Any affine monoid with an active group of invertible elements is isomorphic to a root monoid corresponding to a cone $\sigma$, a $k$-dimensional regular face $\tau$ of $\sigma$, and a collection of Demazure roots
$$
\{e_1^{(r)},e_2^{(r)} \mid r=1,\ldots,k\}
$$
compatible with $\tau$.

Moreover, the group of invertible elements of the root monoid is active if and only if the differences
$$
e_2^{(r)}-e_1^{(r)},
\qquad
r=1,\ldots,k,
$$
are linearly independent.
\end{theorem}

From now on, we fix a root monoid $X_{\sigma,E}$ and assume that the group $G_{\overline{\chi}}$ is active. We now describe the automorphism group of $G_{\overline{\chi}}$.

The center of $G_{\overline{\chi}}$ is
$$
Z(G_{\overline{\chi}})
=
0\times \bigcap_{i=1}^k \Ker \chi_i,
$$
and
$$
\Inn(G_{\overline{\chi}})
\simeq
\BG_a^k\rtimes T',
$$
where
$$
T'
=
T\big/\bigcap_{i=1}^k \Ker \chi_i.
$$

Let $\Psi$ be an automorphism of $G_{\overline{\chi}}$. Since $\BG_a^k$ is the unipotent radical of $G_{\overline{\chi}}$ and all maximal tori in $G_{\overline{\chi}}$ are conjugate, composing $\Psi$ with an inner automorphism we may assume that
$$
\Psi(\BG_a^k)=\BG_a^k
\qquad
\text{and}
\qquad
\Psi(T)=T.
$$
Then $\Psi$ is determined by a pair
$$
(\varphi,\psi)
\in
\Aut(\BG_a^k)\times \Aut(T)
\simeq
\GL_k(\BK)\times \GL_{n-k}(\BZ).
$$

Take arbitrary elements
$$
\alpha,\alpha'\in \BG_a^k,
\qquad
t,t'\in T.
$$
On the one hand,
$$
\Psi\bigl((\alpha,t)\cdot (\alpha',t')\bigr)
=
\Psi\bigl(\alpha+\overline{\chi}(t)\alpha',tt'\bigr)
=
\bigl(
\varphi(\alpha)+
\varphi(\overline{\chi}(t)\alpha'),
\psi(tt')
\bigr).
$$
On the other hand,
$$
\Psi(\alpha,t)\Psi(\alpha',t')
=
(\varphi(\alpha),\psi(t))
(\varphi(\alpha'),\psi(t'))
=
\bigl(
\varphi(\alpha)+
\overline{\chi}(\psi(t))\varphi(\alpha'),
\psi(tt')
\bigr).
$$
Therefore,
$$
\varphi(\overline{\chi}(t)\alpha')
=
\overline{\chi}(\psi(t))\varphi(\alpha')
$$
for all $t\in T$ and $\alpha'\in \BG_a^k$.

Since the characters $\chi_1,\ldots,\chi_k$ are linearly independent, the set of matrices $
\overline{\chi}(t),\ 
t\in T,
$
coincides with the subgroup of all diagonal matrices in $\GL_k(\BK)$. It follows that the operator $\varphi$ is itself diagonal and that $\psi$ preserves every character in the collection $\overline{\chi}$.

Denote by
$$
\Aut(T,\overline{\chi})
$$
the subgroup of automorphisms of $T$ preserving all characters in $\overline{\chi}$. For a diagonal automorphism $\varphi\in \GL_k(\BK)$ and $\psi\in \Aut(T,\overline{\chi})$, the corresponding automorphism $\Psi$ is inner if and only if $\psi$ is the identity automorphism.

We therefore obtain the following proposition.

\begin{proposition}
For an active group $G_{\overline{\chi}}$, the automorphism group of $G_{\overline{\chi}}$ is
$$
\Aut(G_{\overline{\chi}})
=
\Inn(G_{\overline{\chi}})
\rtimes
\Aut(T,\overline{\chi})
$$
and hence
$$
\Aut(G_{\overline{\chi}})
\simeq
\left(
\BG_a^k
\rtimes
\Bigl(
T\big/\bigcap_{i=1}^k \Ker \chi_i
\Bigr)
\right)
\rtimes
\Aut(T,\overline{\chi}).
$$
\end{proposition}

The group $
\Aut(T,\overline{\chi})
$
acts naturally on the lattice $M(\tau)$. Under the induced action of $\Aut(T,\overline{\chi})$ on
$
\BK[G_{\overline{\chi}}]\simeq \BK[X_\tau],
$ and
the characters
$
\chi^{e_1^{(1)}},\ldots,\chi^{e_1^{(k)}}
$
are fixed.

The elements
$$
\chi^{e_1^{(1)}},
\ldots,
\chi^{e_1^{(k)}},
\chi^{q_{k+1}},
\ldots,
\chi^{q_n}
$$
form a basis of the character lattice $M$. Therefore every automorphism
$$
\psi\in \Aut(T,\overline{\chi})
$$
induces an automorphism
$$
\overline{\psi}:M\to M.
$$
The automorphism $\overline{\psi}$ preserves the sublattice $M(\tau)$, fixes all characters in $\overline{\chi}$, and fixes the Demazure roots $
e_1^{(1)},\ldots,e_1^{(k)}.
$
Since
$
\chi_i=e_2^{(i)}-e_1^{(i)},
$
it also fixes the roots
$
e_2^{(1)},\ldots,e_2^{(k)}.
$

The algebra
$$
\BK[X_{\sigma,E}]
\subseteq
\BK[X_\tau]
$$
is invariant under $\psi$ if and only if the induced automorphism $\overline{\psi}$ of $M$ preserves the cone~$\sigma^\vee$. Hence $\psi$ extends to an automorphism of $X_{\sigma,E}$ if and only if $\overline{\psi}$ preserves $\sigma^\vee$.

Denote by
$$
\Aut(M,\sigma,\tau,E)
$$
the group of automorphisms of the lattice $M$ preserving the sublattice $M(\tau)$, the cone $\sigma^\vee$, and the elements of the set $E$. We thus obtain the following theorem.

\begin{theorem}
Let $\sigma\subseteq N_\BQ$ be a strongly convex rational polyhedral cone, let $\tau$ be a regular face of $\sigma$ of dimension $k$, and let
$$
E = \{e_1^{(r)},e_2^{(r)} \mid r=1,\ldots,k\}
$$
be a collection of Demazure roots compatible with $\tau$. Suppose that the characters
$$
\overline{\chi}
=
\{
\chi_i=e_2^{(i)}-e_1^{(i)}
\mid
i=1,\ldots,k
\}
$$
are linearly independent. Then the automorphism group of the corresponding root monoid $X_{\sigma,E}$ is
$$
\Aut(X_{\sigma,E})
=
\Inn(X_{\sigma,E})
\rtimes
\Aut(M, \sigma, \tau,E),
$$
and hence
$$
\Aut(X_{\sigma,E})
\simeq
\left(
\BG_a^k
\rtimes
T/Z(G_{\overline{\chi}})
\right)
\rtimes
\Aut(M, \sigma, \tau,E).
$$
\end{theorem}

\end{document}